\documentclass[twocolumn]{autart}    
\usepackage{graphicx}          
\usepackage{color}

\begin{document}

\begin{frontmatter}
\title{On the regional gradient observability of time fractional diffusion processes}
\thanks[footnoteinfo]{This paper was not presented at any IFAC
meeting. Corresponding author: Y.~Chen, T/F. +1(209)228-4672/4047.}

\author[CUG1,CUG2]{Fudong Ge}\ead{gefd2011@gmail.com},    
\author[UCMerced]{YangQuan Chen}$^{,\dag}$
\ead{yqchen@ieee.org},   
\author[DHU2]{Chunhai Kou}\ead{kouchunhai@dhu.edu.cn}  

\address[CUG1]{School of Computer Science, China University of Geosciences, Wuhan 430074, PR China}
\address[CUG2]{Hubei Key Laboratory of Intelligent Geo-Information Processing, China University of Geosciences, Wuhan 430074, PR China} 
\address[UCMerced]{Mechatronics, Embedded Systems and Automation Lab,
 University of California, Merced, CA 95343, USA }             
\address[DHU2]{Department of Applied Mathematics,
 Donghua University, Shanghai 201620, PR China}        


\begin{abstract}                          
This paper for the first time addresses the concepts of regional gradient observability for the Riemann-Liouville
time fractional order diffusion system in an interested subregion of the whole domain without the knowledge
of the initial vector and its gradient. The Riemann-Liouville time fractional order diffusion system which
replaces the first order time derivative of normal diffusion system by a Riemann-Liouville time fractional
order derivative of order $\alpha\in (0,1]$ is used to well characterize those anomalous sub-diffusion
processes. The characterizations of the strategic sensors when the system under consideration is regional
gradient observability are explored. We then describe an approach leading to the reconstruction of the
initial gradient in the considered subregion with zero residual gradient vector. At last, to illustrate the
effectiveness of our results, we present several application examples where the sensors are zone, pointwise
 or filament ones.
\end{abstract}

\end{frontmatter}

\section{Introduction}
It is well known that the anomalous diffusion processes in various real-world complex systems can be well characterized by using fractional order anomalous
diffusion models (\cite{FCAAMainardi}; \cite{ralf}) after the introduction of continuous time random walks (CTRWs) in \\\cite{CTRWfirst}.
Regarded as a
natural extension of the Brownian motions, the CTRWs are proven to be useful in deriving the time or space fractional order diffusion system  by allowing the
incorporation of waiting time probability density function (PDF) and general jump PDF (\cite{3M}; \cite{Mainardi0}; \cite{Gradenigo}). For example, if the particles
are supposed to jump at fixed time intervals with a incorporating waiting times,  the particles then undergo a sub-diffusion process and the time fractional
 diffusion system is introduced to efficiently describe this process.

As stated in \cite{AEIJAI} and \\ \cite{Ge2016IJC}, instead of analyzing a system by purely theoretical viewpoint (for example, see \cite{curtain}), using the notions of sensors
and actuators to investigate the structures and properties of systems can allow us to understand the system better and consequently enable us to steer the
real-world system in a better way. This situation happens in many real dynamic systems, for example the optimal control
of pest spreading (\cite{caojx}), the flow through porous media microscopic process (\cite{uchaikin2012fractional}), or the swarm of robots moving  through dense forest (\cite{2spears}) etc. It is now widely believed that fractional order controls can offer better performance not achievable before
using integer order controls systems (\cite{Mandelbrot}; \cite{Torvik}). This is the reason why the fractional order models are superior in comparison with the integer
order models. Moreover, it is worth noting that in many real dynamic systems, the regional observation problem occurs naturally  when one is interested in the
knowledge of the states in a subregion of the spatial domain (\cite{AEIJAI}; \cite{flux1}; \cite{rgb2}). Focusing on regional observations would allow for a reduction in the number of physical sensors and offer the potential to reduce computational requirements in some cases. In addition, it should be pointed out that the concepts of regional observability are of great help to reconstruct the initial vector for those non-observable system when we are interested in the knowledge of the initial vector only in a critical subregion of the system domain.

\begin{figure}
\begin{minipage}[t]{1\linewidth}
\centering
\includegraphics[width=0.8\textwidth]{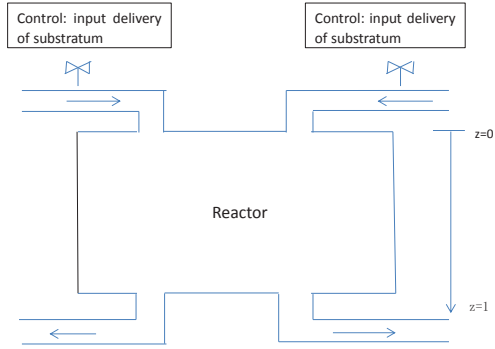}
\end{minipage}%
\caption{Regulation of the concentration flux of the substratum at the upper bottom of the reactor}
\label{fig1}
\end{figure}

Motivated by the argument above, in this paper, by considering the locations, number and spatial distributions of sensors, our goal is to study the regional
gradient observability of the Riemann-Liouville time fractional order diffusion process, which is introduced to better characterize  those
sub-diffusion processes (\cite{henry2000fractional}.)
More precisely, consider the problem $(\ref{problem})$ below and suppose that the initial
vector $y_0$ and its gradient $\nabla y_0$ are unknown and the measurements are given by using output functions (depending on
the number and structure of sensors). The purpose here is to reconstruct the initial gradient vector $\nabla y_0$ on a given subregion of
the whole domain of interest. We also explore the characterizations of strategic sensors when the system is regional gradient observability.
Moreover, there are many applications of gradient modeling.  For example, the concentration regulation
of a substrate at the upper bottom of a  biological reactor sub-diffusion process,
which is observed between two levels (See Fig. \ref{fig1});
Anther example the energy exchange problem  between a casting plasma on a plane target which is perpendicular to the direction of the flow sub-diffusion process from measurements carried out by internal thermocouples (\cite{rgb1}).
For richer background on gradient modeling, we refer the reader to \cite{Jorge} and \cite{gradientmodel}.  To the best
of our knowledge, no results are available on this topic and we hope that the results obtained here
could provide some insights into the control theory analysis of the fractional order diffusion systems and be useful in real-life applications.

The rest contents of the present paper are structured  as follows. The problem studied and some preliminaries are introduced in the next section and in  Section~$\ref{section3},$
we focus on the characteristic of the strategic sensors. An  approach which enables us to reconstruct the initial gradient vector of the system under consideration
 in the considered subregion is addressed in Section~$\ref{section4}$. Several application examples are worked out in the end for illustrations.

\section{Problem formulation and preliminaries}

In this section, we formulate the regional gradient observability problems for the  Riemann-Liouville time fractional order diffusion system and then
introduce some preliminary results to be used thereafter.
\subsection{Problem formulation}
Let $\Omega$ be a connected, open bounded subset of $\mathbf{R}^n$ with Lipschitz continuous boundary $\partial\Omega$ and
consider the following abstract time fractional diffusion process:
\begin{equation}\label{problem}
\left.\begin{array}{l}
_0D^{\alpha}_{t}y(t)=Ay(t), ~~ t\in [0,b],~0<\alpha\leq 1,
\\\lim\limits_{t\to 0^+}{} _0I^{1-\alpha}_{t}y(t)=y_0\mbox{ supposed to be unknown,}
\end{array}\right\}
\end{equation}
where $A$ generates a strongly continuous semigroup $\{\Phi(t)\}_{t\geq 0}$ on the Hilbert space $Y: =L^2(\Omega)$, $-A$ is a uniformly elliptic operator, $y\in  L^2(0,b; Y)$,
$_0D^{\alpha}_t$ and ${}_0I^{\alpha}_t $ denote the Riemann-Liouville fractional order derivative and integral with respect to time $t$, respectively,
given by (\cite{Kilbas} and \cite{Podlubny})
\begin{eqnarray*}
_0D_t^\alpha y(t)&=&\frac{\partial}{\partial t}{}_0I_t^{1-\alpha }y(t),~0< \alpha\leq 1~~ \mbox{ and }\\
_0I^{\alpha}_ty(t)&=&\frac{1}{\Gamma(\alpha)}\int^{t}_{0}{(t-s)^{\alpha-1}y(s)ds},~ \alpha>0.
\end{eqnarray*}
The measurements (possibly unbounded) are given depending on the number and the structure of the sensors with dense domain in $L^2(0,b;Y)$ and
range in $L^2(0,b;\mathbf{R}^p)$ as follows:
\begin{eqnarray}\label{outputfunction}
z(t)=Cy(t),
\end{eqnarray}
where $p\in \mathbf{N}$ is the finite number of sensors.

Let  $y_0\in H_0^1(\Omega)$   and both the initial vector  $y_0$ and its gradient are supposed to be unknown.
The system $(\ref{problem})$ admits a unique mild solution given by  (\cite{Ge2016JMAA} and \cite{liu2015siam}):
\begin{eqnarray}\label{solution}
y(t)=S_\alpha(t)y_0,~~t\in [0,b],
\end{eqnarray}
where $S_\alpha(t)=\alpha t^{\alpha-1}\int_0^\infty{\theta\phi_\alpha(\theta)\Phi(t^\alpha\theta)}d\theta,$
 $\{\Phi(t)\}_{t\geq 0}$ is the strongly continuous semigroup generated by $A$, $\phi_\alpha(\theta)= \frac{1}{\alpha}\theta
^{-1-\frac{1}{\alpha}}\psi_\alpha(\theta^{-\frac{1}{\alpha}})$ and  $\psi_\alpha$ is a probability density function defined by $(\theta\in (0,\infty))$
\begin{eqnarray}\label{2.7}
\psi_\alpha(\theta)=\frac{1}{\pi}\sum\limits_{n = 1}^\infty {( - 1)^{n - 1}\theta^{-\alpha n-1}\frac{\Gamma(n\alpha+1)}{n!}
\sin(n\pi \alpha)}\end{eqnarray}
satisfying (\cite{Mainardi})
\begin{eqnarray*}\int_0^\infty{\psi_\alpha(\theta)}d\theta=1 \mbox{ and } \int_0^\infty{\theta^\nu
\phi_\alpha(\theta)}d\theta=\frac{\Gamma(1+\nu)}{\Gamma(1+\alpha\nu)}, \nu\geq 0.\end{eqnarray*}
\subsection{Definitions and characterizations}
Let $\omega \subseteq \Omega$ be a given region of positive Lebesgue measure and let
\begin{equation}
y_0=\left\{\begin{array}{l}
y_0^1,~\omega~\mbox{ to be estimated,}
\\ y_0^2,~\Omega\backslash \omega~ \mbox{ undesired.}
\end{array}\right.
\end{equation}
Then the regional gradient observability problem
is concerned with the directly reconstruction of the initial gradient vector $\nabla y_0^1$ in $\omega$.
Consider the following two restriction mappings
\begin{eqnarray*}
p_\omega:~\left\{\begin{array}{l}
\left(L^2(\Omega)\right)^n\to \left(L^2(\omega)\right)^n_,
\\~ \xi\to \xi|_{\omega}
\end{array}\right.
p_{1\omega}:~\left\{\begin{array}{l}
L^2(\Omega)\to L^2(\omega),
\\~ y\to y|_{\omega}.
\end{array}\right.
\end{eqnarray*}

Their adjoint operators  are,  respectively,  denoted by
\begin{equation}
p_\omega^*:~\left\{\begin{array}{l}
\left(L^2(\omega)\right)^n \to \left(L^2(\Omega)\right)^n,
\\~ \xi\to p_\omega^* \xi=\left\{\begin{array}{l}\xi , ~~x\in \omega ,  \\0,~~x\in \Omega\backslash \omega\end{array}\right.
\end{array}\right.
\end{equation}
and
\begin{equation}
p_{1\omega}^*:~\left\{\begin{array}{l}
L^2(\omega)\to L^2(\Omega),
\\~ y\to p_{1\omega}^* y=\left\{\begin{array}{l}y , ~~x\in \omega ,  \\0,~~x\in \Omega\backslash \omega.\end{array}\right.
\end{array}\right.
\end{equation}
Moreover, by Eq. $(\ref{solution})$, the output function  $(\ref{outputfunction})$ gives
\begin{eqnarray}
z(t)=CS_\alpha(t)y_0=K(t)y_0,
\end{eqnarray}
where $K:H_0^1(\Omega) \to L^2(0,b;\mathbf{R}^p)$. To obtain the adjoint operator of $K$, we have\\
\textbf{Case 1. $C$ is bounded~(e.g. zone sensors) }\\
Denote the adjoint operator of $C$ and $S_\alpha$ by $C^*$ and $S_\alpha^*$, respectively.
Since $S_\alpha$ is a bounded operator (\cite{zhouyong}), we get that the adjoint operator of $K$ can be given by
\begin{equation}
K^*:~\left\{\begin{array}{l}L^2(0,b;\mathbf{R}^p)\to Y,
\\~ z\to \int_0^b{S^*_\alpha(s)C^*z(s)}ds.
\end{array}\right.
\end{equation}
\textbf{Case 2. $C$ is unbounded~(e.g. pointwise sensors) }\\
Note that $C$ is densely defined, then $C^*$ exists. To state our results, the following two assumptions are needed:

$(A_1)$  $CS_\alpha(t)$ can be extended to  a bounded linear operator $\overline{CS_\alpha(t)}$ in $\mathcal{L}\left(Y,L^2(0,b;\mathbf{R}^p)\right)$;\\
$(A_2)$  $(CS_\alpha)^*$ exists and $(CS_\alpha)^*=S_\alpha^*C^*$.

Extend $K$ by $K(t) y_0= \overline{CS_\alpha(t)}y_0, $  one has $K\in \mathcal{L}\left(Y,L^2(0,b;\mathbf{R}^p)\right).$
Based on the Hahn-Banach theorem,  similar to the argument in \cite{dualityrelationship}, it is possible to derive the duality theorems
as in \cite{curtain} and \cite{dolecki} with the above two assumptions.
Then the adjoint operator of $K$ can be defined as
\begin{equation}
K^*:~\left\{\begin{array}{l}
D(K^*)\subseteq L^2(0,b;\mathbf{R}^p)\to Y,
\\~ z\to \int_0^b{S^*_\alpha(s)C^*z(s)}ds.
\end{array}\right.
\end{equation}
Let $\nabla: H_0^1(\Omega) \to \left(L^2(\Omega)\right)^n $ be an operator defined by
\begin{eqnarray}
y \to \nabla y(x):=\left(\frac{\partial y}{\partial x_1},\frac{\partial y}{\partial x_2},\cdots,\frac{\partial y}{\partial x_n}\right).
\end{eqnarray}
We see that the adjoint of the gradient operating on a connected, open  bounded subset $\Omega$ with a Lipschitz continuous
boundary $\partial\Omega$ is minus the divergence operator, i.e.,  $\nabla^*: \left(L^2(\Omega)\right)^n \to H^{-1}(\Omega)$
is given by (\cite{duality})
\begin{eqnarray}\label{2.10}
\xi \to \nabla^* \xi:=v,
\end{eqnarray}
where $v$ solves the following Dirichlet problem
\begin{equation}
\left\{\begin{array}{l}
v=- \mbox{div}(\xi) ~\mbox{ in } \Omega,
\\~v=0 ~\mbox{ on } \partial\Omega.
\end{array}\right.
\end{equation}

%

Similar to the discussion in \cite{curtain}; \cite{dolecki} and \cite{dualityrelationship}, it follows that the
necessary and sufficient condition for the regional weak observability of the system described by
$(\ref{problem})$ and $(\ref{outputfunction})$ in $\omega$ at time $b$
is that
$
Ker\left( Kp_{1\omega}^*\right)=\{0\}\subseteq  L^2(\Omega)
$
and we see the following definition.
\begin{defn}\label{rgodef}
The system $(\ref{problem})$ with output function $(\ref{outputfunction})$ is said to be
 regional weak gradient observability in $\omega$ at time $b$ if and only if
\begin{eqnarray}
Ker\left( K\nabla^*p_{\omega}^*\right)=\{(0,0,\cdots,0)\}\subseteq \left(L^2(\Omega)\right)^n.
\end{eqnarray}
\end{defn}

\begin{prop}\label{proposition1}  There is an equivalence among the following properties:

$\left<1\right>$ The system $(\ref{problem})$ is regional weak  gradient observability
 in $\omega$ at time $b$;
\\
$\left<2\right>$ $\overline{Im \left(p_{\omega} \nabla K^*\right)}= \left(L^2(\omega)\right)^n;$
\\
$\left<3\right>$
The operator $p_\omega \nabla K^* K \nabla^*p^*_\omega$ is positive definite.
\end{prop}
\textbf{Proof.}
By Definition $\ref{rgodef}$, it is obvious to know that  $\left<1\right>\Leftrightarrow \left<2\right>.$
%
%
%
As for  $\left<2\right>\Leftrightarrow \left<3\right>,$ in fact, we have
\begin{eqnarray*}
\overline{Im\left( p_{\omega}\nabla K^*\right)}=\left(L^2(\omega)\right)^n
&\Leftrightarrow&
(p_\omega \nabla K^*z,y)=0,~\forall z\in Z\Rightarrow\\&~~& y=(0,0,\cdots,0)\in \left(L^2(\Omega)\right)^n.
\end{eqnarray*}
Let $z=K\nabla^*p_\omega^* y \in  Z$, which then allows us to complete the  proof.

\begin{rem}
$\\$
$\left.1\right)$ When $\alpha=1$, the system $(\ref{problem})$ is deduced to the normal diffusion process as considered in \cite{rgb2},
which is a particular case of our results.

$\left.2\right)$ A system which is gradient observable on $\omega$ is gradient observable on $\omega_1$ for every $\omega_1\subseteq \omega.$
Moreover, the Definition $\ref{rgodef}$ is also valid for the case when $\omega=\Omega$ and  there
exist systems that are not gradient observable but regionally  gradient observable. This can be illustrated by the
following example.
\end{rem}

\subsection{An example}
Let $\Omega=[0,1]\times [0,1] \subseteq \mathbf{R}^2$ and consider the following time fractional order diffusion system of order $\alpha\in (0, 1]$.
\begin{eqnarray}\label{counterexample}
\left\{\begin{array}{l}
_0D^{\alpha}_{t}y(x_1,x_2,t)=\left(\frac{\partial^2}{\partial x_1^2}+\frac{\partial^2}{\partial x_2^2}\right)y(x_1,x_2,t)\mbox{ in }\Omega \times [0,b],
\\ y(\xi,\eta, t)= 0~\mbox{ on }\partial\Omega \times [0,b],\\
\lim\limits_{t\to 0^+}{}_0I^{1-\alpha}_{t}y(x_1,x_2,t)=y_0(x_1,x_2)~\mbox{ in }\Omega
\end{array}\right.
\end{eqnarray}
with the output functions
\begin{eqnarray}
z(t)=Cy(t)=\int_0^1\int_0^1{f(x_1,x_2)y(x_1,x_2,t)}dx_1dx_2,
\end{eqnarray}
where $f(x_1,x_2)=\delta(x_1-1/2) \sin(\pi x_2)$ and $\delta(x)$ is the Dirac delta function on the real number line that is zero everywhere except at zero.

According to  the problem $(\ref{problem})$, $A=\frac{\partial^2 }{\partial x_1^2}+\frac{\partial^2 }{\partial x_2^2}$. Then the eigenvalue,
 eigenvector and the semigroup $\Phi(t)$ on $Y$ generated by $A$ are respectively
 $\lambda_{ij}=-(i^2+j^2)\pi^2$,  $\xi_{ij}(x_1,x_2)=2\sin (i\pi x_1)\sin (j\pi x_2)$
and
$\Phi(t)y(x)=\sum\limits_{i,j=1}^{\infty}{\exp(\lambda_{ij}t)(y,\xi_{ij})\xi_{ij}}(x).$
Moreover,  one has (\cite{Ge2016JMAA})
\[S_\alpha(t)y_0(x)=\sum\limits_{i,j=1}^{\infty}t^{\alpha-1}E_{\alpha,\alpha}(\lambda_{ij}t^{\alpha})(y_0,\xi_{ij})\xi_{ij}(x),
\]
where $E_{\alpha,\beta}(z):=\sum\limits_{i=0}^\infty
{\frac{z^i}{\Gamma(\alpha i+\beta)}},$ $\mathbf{Re}{\kern 2pt}\alpha>0, ~\beta,z\in \mathbf{C}$
 is the generalized Mittag-Leffler function in two parameters.

Next, we show that there is a gradient vector $g$, which is not gradient observable in
the whole domain but gradient observable in a subregion $\omega\subseteq \Omega$.

Let  $ g=\frac{1}{\pi}(\cos(\pi x_1)\sin(3 \pi x_2),3 \sin(\pi x_1)\cos(3 \pi x_2))\in \left(L^2(\Omega)\right)^2$. By Eq. $(\ref{2.10})$,
we obtain that $\nabla^*g=10\sin(\pi x_1)\sin(3 \pi x_2)$,
then
\begin{eqnarray*}
K\nabla^*g&=& CS_\alpha(t)\nabla^*g\\
&=& 40\sum\limits_{i,j=1}^{\infty}t^{\alpha-1}
E_{\alpha,\alpha}(\lambda_{ij}t^{\alpha})\int_0^1{\sin(\pi x_1)\sin (i\pi x_1)}dx_1
\\&~&\times\int_0^1{\sin(3 \pi x_2)\sin (j\pi x_2)}dx_2\\
&~&\times \sin \left(\frac{i\pi}{2}\right)\int_0^1{\sin (\pi x_2)\sin (j\pi x_2)}dx_2\\
&=&0.
\end{eqnarray*}
However, let $\omega= [0,1]\times[0,1/6],$ we see that
\begin{eqnarray*}
K\nabla^*p^*_\omega p_\omega g&=& CS_\alpha(t)\nabla^*p^*_\omega p_\omega g\\
&=&40 \sum\limits_{i,j=1}^{\infty}t^{\alpha-1}
E_{\alpha,\alpha}(\lambda_{ij}t^{\alpha})\int_0^1{\sin(\pi x_1)\sin (i\pi x_1)}dx_1
\\&~&\times\int_{0}^{1/6}{\sin(3 \pi x_2)\sin (j\pi x_2)}dx_2\\
&~&\times \sin \left(\frac{i\pi}{2}\right)\int_0^1{\sin (\pi x_2)\sin (j\pi x_2)}dx_2\\
&=& \frac{5\sqrt{3}t^{\alpha-1}E_{\alpha,\alpha}(-2\pi^2t^{\alpha})}{8\pi }\neq 0,
\end{eqnarray*}
which means that $g $ is gradient observable in $\omega.$

The following two lemmas play a significant role to obtain our results.
\begin{lem}(\cite{podlubnyadjoint})\label{integral by part}
For any $t\in [a,b]$,  $\alpha\in (0,1)$, the following formula holds
 \begin{eqnarray*}\label{derivativeparts}
\int_a^b{f(t){}_aD_t^{\alpha}g(t)}dt=
\left[f(t){}_aI_t^{1-\alpha}g(t)\right]_{t=a}^{t=b}-
\int_a^b{g(t){}_t^CD_b^{\alpha}f(t)}dt.
 \end{eqnarray*}
\end{lem}
\begin{lem}(\cite{Bernard})\label{Lem1}
Let $\Omega \subseteq \mathbf{R}^n$ be
an open set and $C_0^{\infty}(\Omega)$ be the class of  infinitely differentiable functions on $\Omega$
 with compact support in $\Omega$ and $u\in L^1_{loc}(\Omega)$ be such that
\begin{eqnarray}
~~~~~~~~~~~~~~\int_\Omega u(x)\psi(x)dx=0,~~~\forall\psi \in C_0^{\infty}(\Omega).
\end{eqnarray}
Then $u=0$ almost everywhere in $\Omega.$
\end{lem}

\section{Regional strategic sensors}
\label{section3}
This section is devoted to addressing the characteristic of sensors when the studied system is regionally gradient
observable in a given subregion of the whole domain.

Firstly,  we recall that a sensor can be defined by a couple $(D,f)$ where $D\subseteq \Omega $ is the support
of the sensor and $f$ is its spatial distribution. For example, if  $D=\{\sigma\}$ with  $\sigma\in \bar{\Omega}$ and $f=\delta_\sigma,$ where $\delta_\sigma=\delta(\cdot-\sigma)$ is
the Dirac delta function in $\omega$ at time $b$ that is zero everywhere except at $\sigma,$
the sensor is called pointwise sensor. In this case the operator $C$ is unbounded and the output function can be written as
$
z(t)=y(\sigma,t).
$
It is called zone sensor when $D\subseteq \bar{\Omega}$ and $f\in Y$. The output function is bounded and can be defined as follows:
$
z(t)=\int_D{y(x,t)f(x)}dx.
$
For more information on the structure characteristic and properties of sensors and actuators, we refer the reader to
(\cite{AEIJAI2}, \cite{AEIJAI}, \cite{rgb2}) and the references cited therein.

Next, to state our results,  it is supposed that the measurements are made by $p$ sensors $(D_i,f_i)_{1\leq i\leq p}$,
where $D_i\subseteq \Omega$ and $f_i\in L^2(\Omega)$, $i=1,2,\cdots,p.$
Then  $(\ref{problem})$ can be rewritten as
\begin{equation}\label{problem22}
\left\{\begin{array}{l}
{}_0D^{\alpha}_{t}y(x,t)=Ay(x,t)~\mbox{ in }\Omega \times [0,b],\\
y(\eta,t)=0~~\mbox{ on }  \partial \Omega \times [0,b],\\
\lim\limits_{t\to 0^+}{} _0I^{1-\alpha}_{t}y(x,t)=y_0(x)~\mbox{ in }\Omega
\end{array}\right.
\end{equation}
with the measurements
\begin{eqnarray}
z(t)=Cy(x,t)=\left(z_1(t),z_2(t),\cdots,z_p(t)\right)^T\in \mathbf{R}^p,
\end{eqnarray}
where $z_i(t)=\left(y(\cdot,t),f_i\right)_{L^2(D_i)}$.

Moreover, since the operator  $-A$ is a uniformly elliptic operator, for any $y_i\in L^2(0,b;Y),$ $i=1,2$, $A$ satisfies
\begin{eqnarray*}
\int_Q{y_1(x,t)Ay_2(x,t)}dtdx-\int_Q{y_2(x,t) A^*y_1(x,t)}dtdx\\
=\int_{\partial\Omega \times [0,b]}{\left[y_1(\eta,t)\frac{\partial y_2(\eta,t)}{\partial v_A}-y_2(\eta,t) \frac{\partial y_2(\eta,t)}{\partial v_{A^*}}\right]}dtd\eta,
\end{eqnarray*}
where $A^*$ is the adjoint operator of $A$. Moreover, by \cite{Hilbert}, there exists a sequence $(\lambda_j,\xi_{jk}): k=1,2,\cdots,r_j$, $j=1,\cdots $ such that

$1)$ Each $\lambda_j$ $( j=1,2,\cdots)$ is the eigenvalue of the operator $A$ with multiplicities $r_j$ and
\[0> \lambda_1>\lambda_2> \cdots> \lambda_j>\cdots, ~~\lim\limits_{j\to \infty}\lambda_j=-\infty.\]

$2)$ For each $ j=1,2,\cdots$,  $\xi_{jk}~(k=1,2,\cdots,r_j)$ is the orthonormal eigenfunction corresponding to $\lambda_j$, i.e.,
\[
(\xi_{jk_m},\xi_{jk_n})=\left\{\begin{array}{l}1,~~k_m=k_n,
\\
0,~~k_m\neq k_n,
\end{array}\right.\]
where $1\leq k_m,k_n\leq  r_j,$ $k_m,k_n\in \mathbf{N}$ and  $(\cdot,\cdot)$ is the inner product of space $Y$.
Then it follows that the strongly continuous semigroup $\{\Phi(t)\}_{t\geq 0}$ on  $Y$ generated by $A$ can be expressed as
\begin{eqnarray}\label{phi}
\Phi(t)y(x)=\sum\limits_{j=1}^{\infty}{\sum\limits_{k=1}^{r_j}{\exp(\lambda_jt)(y,\xi_{jk})\xi_{jk}(x)}}, ~~x\in \Omega,
\end{eqnarray}
the sequence $\{\xi_{jk}, k=1,2,\cdots,r_j, j=1,2,\cdots \}$ is an  orthonormal basis in $Y$ and   for any $y_*(x)\in Y$, it can be expressed as
$y_*(x)=\sum\limits_{j=1}^{\infty}{\sum\limits_{k=1}^{r_j}(y_*,\xi_{jk})\xi_{jk}(x)}.$

\begin{defn}
A sensor (or a suite of sensors) is said to be gradient $\omega-$strategic if the observed system is regionally gradient observable in $\omega$.
\end{defn}
\begin{lem}\label{lemma3.1}
For any $z(t)=\left(z_1(t),z_2(t),\cdots,z_p(t)\right)^T\in \mathbf{R}^p$ with $z_i\in L^2(0,b),$ $i=1,2,\cdots,p$, suppose that
 $e(x,t)$ satisfies the following system
\begin{equation}\label{efunction}
\left\{\begin{array}{l}
{}^C_tD^{\alpha}_{b}e(x,t)=-A^*e(x,t)+\sum\limits_{i=1}^p{p_{D_i}f_i(x)z_i(t)} \\
{\kern 150 pt}\mbox{ in } \Omega \times [0,b],\\
e(\eta,t)=0~~\mbox{ on }  \partial \Omega \times [0,b],
\\ e(x,b)=0 ~~\mbox{ in } \Omega,
\end{array}\right.
\end{equation}
where $A^*$ is the adjoint operator of $A$ and $_t^CD^{\alpha}_b$ denotes the right-sided Caputo fractional order derivative with respect
to time $t$ of order $\alpha\in (0,1]$ given by (\cite{Kilbas}; \cite{Podlubny} and \cite{podlubnyadjoint})
\begin{eqnarray}
_t^CD_b^\alpha e(x,t)=\frac{-1}{\Gamma(1-\alpha)}\int_t^b{(\tau-t)^{-\alpha}\frac{\partial}{\partial \tau}e(x,\tau)}d\tau.
\end{eqnarray}
Then we obtain that
$
K^*z=-e(x,0).
$
\end{lem}
\textbf{Proof.}
Replacing $t$ in $_t^CD_b^\alpha e(x,t)$ by $b-t$, we have
\begin{eqnarray*}
\begin{array}{l}
~~~~{}^C_{(b-t)}D^{\alpha}_{b}e(x,b-t)\\=\frac{-1}{\Gamma(1-\alpha)}\int_{b-t}^b{(\tau-b+t)^{-\alpha}\frac{\partial}{\partial \tau} e(x,\tau)}d\tau\\
=\frac{-1}{\Gamma(1-\alpha)}\int_0^t{(t-s)^{-\alpha}\left[-\frac{\partial}{\partial s} e(x,b-s)\right]}ds\\
=-{}^C_0D^{\alpha}_{t}e(x,b-t).
\end{array}\end{eqnarray*}
Then  the system  $(\ref{efunction})$ is equivalent to
\begin{equation}
\left\{\begin{array}{l}
{}^C_0D^{\alpha}_te(x,b-t)=A^*e(x,b-t) ~ \\
{\kern 41pt}-\sum\limits_{i=1}^p{p_{D_i}f_i(x)z_i(b-t)}\mbox{ in } \Omega \times [0,b],\\
e(\eta,b-t)=0~~\mbox{ on }  \partial \Omega \times [0,b],
\\ e(x,b-0)=0 ~~\mbox{ in } \Omega.
\end{array}\right.
\end{equation}
Similar to  the argument in \cite{Ge2016JMAA} and  \cite{Japan}, the mild solution of  $(\ref{efunction})$ can be given by
\begin{eqnarray*}
e(x,t)=-\int_0^{b-t}{S_\alpha^*(b-t-s)\sum\limits_{i=1}^p{p_{D_i}f_i(x)z_i(b-s)} }ds.
\end{eqnarray*}
On the other hand, it follows from the definition of the adjoint operator of $K$ that
\begin{eqnarray*}
K^*z=\int_0^b{S^*_\alpha(s)C^*z(s)}ds=\int_0^b{S^*_\alpha(s)\sum\limits_{i=1}^p{p_{D_i}f_i(x)z_i(s)}}ds.
\end{eqnarray*}
This allows us to complete the proof.

\begin{thm}\label{theorem3.1} For any $j=1,2,\cdots$, $s=1,2,\cdots,n$, given arbitrary $b>0$, define the following $p\times r_j$ matrices $G_j^s$
\begin{equation} \label{G}
G_j^s=\left[ {\begin{array}{*{20}{c}}
{\xi_{j1}^{1s}}&{\xi_{j2}^{1s}}&{\cdots}&{\xi_{jr_j}^{1s}}\\
{\xi_{j1}^{2s}}&{\xi_{j2}^{2s}}&{\cdots}&{\xi_{jr_j}^{2s}}\\
{\vdots}&{\vdots}&{\vdots}&{\vdots}\\
{\xi_{j1}^{ps}}&{\xi_{j2}^{ps}}&{\cdots}&{\xi_{jr_j}^{ps}}
\end{array}} \right]_{p\times r_j},\end{equation}
where $\xi_{jk}^{is}=\left(\frac{\partial \xi_{jk}}{\partial x_s},f_i\right)_{L^2(D_i)}$, $i=1,2,\dots,p$ and  $k=1,2,\cdots,r_j$.
For all $ j=1,2,\cdots,$ let $y_{jks}=(p_{1\omega}^*y_s,\xi_{jk})$ and $y_{js}=(y_{j1s},y_{j2s},\cdots,y_{jr_js})^T\in \mathbf{R}^{r_j}$.
Then the necessary and sufficient condition for the gradient $\omega-$strategic of the sensors $(D_i,f_i)_{1\leq i\leq p}$ is that
\begin{eqnarray*}
\sum\limits_{s=1}^n G_j^sy_{js}=\mathbf{0}_p:=(0,0,\cdots,0)^T\in \mathbf{R}^p ~\Rightarrow{\kern 60pt}\\
{\kern 80pt}y =\mathbf{0}_n: =(0,\cdots,0)^T\in \left( L^2( \omega)\right)^n.
\end{eqnarray*}
In particular, when $n=1,$ the sensors $(D_i,f_i)_{1\leq i\leq p}$ is  gradient $\omega-$strategic if and only if
\begin{eqnarray*}
&~&(1)~p\geq r=\max\{r_j\}; \\
&~&(2)~rank ~G_j^1=r_j \mbox{ for all } j=1,2,\cdots.\end{eqnarray*}
\end{thm}
$\mathbf{Proof.}$
Given arbitrary $b>0$, by Definition $\ref{rgodef}$,
the sensors $(D_i,f_i)_{1\leq i\leq p}$ are gradient $\omega-$strategic if and only if for any $z\in L^2(0,b;\mathbf{R}^p)$,
\begin{eqnarray*}
\left\{y\in \left(L^2(\omega)\right)^n|   \left(p_\omega\nabla K^*z,y\right)_{\left(L^2(\omega)\right)^n}=0\right\}
 \Rightarrow y=\mathbf{0}_n,
\end{eqnarray*}
where $y=(y_1,y_2,\cdots,y_n)$ with $y_s\in L^2(\omega)$.

Let $x=(x_1,x_2, \cdots,x_n)\in \Omega$. By Lemma $\ref{lemma3.1},$ we have
\begin{eqnarray*}
\left(p_\omega\nabla K^*z,y\right)_{\left(L^2(\omega)\right)^n}
&= & \left(\nabla K^*z,p_\omega^* y\right)_{\left(L^2(\Omega)\right)^n}\\
&=&\sum\limits_{s=1}^n\left(\frac{ \partial(K^*z)}{\partial x_s},p_{1\omega}^* y_s\right)_{L^2(\Omega)}\\
&=&\sum\limits_{s=1}^n\left(\frac{ \partial \left[-e(x,0)\right]}{\partial x_s},p_{1\omega}^* y_s\right)_{L^2(\Omega)},
\end{eqnarray*}
where $e$ is the solution of the system $(\ref{efunction})$.

Next, we explore the exact expression of $\left(p_\omega\nabla K^*z,y\right)_{\left(L^2(\omega)\right)^n}.$
Consider the following problem
\begin{equation}\label{adjointproblem}
\left\{\begin{array}{l}
_0D^{\alpha}_t\rho(x,t)=A\rho(x,t) ~\mbox{in }\Omega\times [0,b],
\\ \rho(\eta,t)=0~~\mbox{on }\partial\Omega\times [0,b],
\\ \lim\limits_{t\to 0^+}{}_0I^{1-\alpha}_t\rho(x,t)=p_{1\omega}^*y_s(x){\kern 5pt}~\mbox{in }\Omega,
\end{array}\right.
\end{equation}
where $s=1,2,\cdots,n$ and the unique mild solution of system $(\ref{adjointproblem})$ can be given by (\cite{Japan})
\begin{eqnarray}
\begin{array}{l}
\rho(x,t)=\sum\limits_{j=1}^{\infty}\sum\limits_{k=1}^{r_j}
 E_{\alpha,\alpha}(\lambda_jt^{\alpha})\left(p_{1\omega}^*y_s,\xi_{jk}\right) \xi_{jk}(x).
\end{array}\end{eqnarray}
 Multiplying both sides of $(\ref{efunction})$ with  $\frac{\partial\rho(x,t)}{\partial x_s}$ and integrating
 the results over the domain $Q:=\Omega\times [0,b]$,
\begin{eqnarray*}
\begin{array}{l}
\int_Q{\left[{}_t^CD^{\alpha}_{b}e(x,t)\right]\frac{\partial\rho(x,t)}{\partial x_s}}dtdx\\
=-\int_Q{A^*e(x,t)\frac{\partial\rho(x,t)}{\partial x_s}}dtdx\\{\kern 8pt}
+\int_0^b \int_\Omega {\sum\limits_{i=1}^p{p_{D_i}f_i(x)z_i(t)} \frac{\partial\rho(x,t)}{\partial x_s} }dxdt .
\end{array}\end{eqnarray*}
Consider the fractional integration by parts $(\ref{derivativeparts})$ in Lemma $\ref{integral by part}$,  one has
\begin{eqnarray*}
\begin{array}{l}
\int_Q{\left[{}_t^CD^{\alpha}_be(x,t)\right] \frac{\partial\rho(x,t)}{\partial x_s}}dtdx\\=
-\int_\Omega{\left[\lim\limits_{t\to 0^+}{}_0I^{1-\alpha}_t\rho (x,t)\right] \frac{\partial e(x,0)}{\partial x_s}}dx\\{\kern 8pt}-\int_Q{e(x,t)\left[ {}_0D^{\alpha}_t\frac{\partial\rho(x,t)}{\partial x_s}\right]}dtdx.
\end{array}\end{eqnarray*}
Then the boundary condition gives
 \begin{eqnarray*}
 \begin{array}{l}
\int_\Omega{ \frac{\partial \left[-e(x,0)\right]}{\partial x_s}\left[\lim\limits_{t\to 0^+}{}_0I^{1-\alpha}_te(x,t)\right]}dx\\\textcolor[rgb]{1.00,0.00,0.00}{=}
\int_0^b \int_\Omega {\sum\limits_{i=1}^p{p_{D_i}f_i(x)z_i(t)} \frac{\partial\rho(x,t)}{\partial x_s} }dxdt ,~s=1,2,\cdots,n.
\end{array}\end{eqnarray*}
Thus, we have
\begin{eqnarray*}
&~&\left(p_\omega\nabla K^*z,y\right)_{\left(L^2(\omega)\right)^n}\\
&=&\sum\limits_{s=1}^n\left(\frac{ \partial \left[-e(x,0)\right]}{\partial x_s},p_{1\omega}^* y_s\right)_{L^2(\Omega)}\\
&=&\sum\limits_{s=1}^n\left(\frac{ \partial \left[-e(x,0)\right]}{\partial x_s},\left[\lim\limits_{t\to 0^+}{}_0I^{1-\alpha}_te(x,t)\right]\right)_{L^2(\Omega)}\\
&=&\sum\limits_{s=1}^n\sum\limits_{j=1}^{\infty}\sum\limits_{k=1}^{r_j}\sum\limits_{i=1}^{p}\int_0^b{
 E_{\alpha,\alpha}(\lambda_j\tau^{\alpha})z_i(\tau)}d\tau \left(\frac{\partial \xi_{jk}}{\partial x_s},p_{D_i}f_i\right) y_{jks}.
\end{eqnarray*}
By Lemma $\ref{Lem1}$, since $z\in L^2(0,b;\mathbf{R}^p)$ is arbitrary, we see that the system $(\ref{problem})$ is regionally gradient observable
in $\omega$ at time $b$ if and only if
\begin{eqnarray}\label{3.12}
\sum\limits_{s=1}^n \sum\limits_{j=1}^{\infty}E_{\alpha,\alpha}(\lambda_jt^{\alpha})\sum\limits_{i=1}^{p}
\sum\limits_{k=1}^{r_j}\xi_{jk}^{is} y_{jks} =\mathbf{0}_p
 \Rightarrow y = \mathbf{0}_n,
\end{eqnarray}
i.e., for any $y=(y_1,y_2,\cdots,y_n)\in \left(L^2(\Omega)\right)^n$, one has
 \begin{eqnarray}\label{3.12}
\sum\limits_{j=1}^{\infty}{E_{\alpha,\alpha}
(\lambda_jt^\alpha)}\sum\limits_{s=1}^n G_j^sy_{js}=\mathbf{0}_p~\Rightarrow y = \mathbf{0}_n,
\end{eqnarray}
where $y_{js}=(y_{j1s},y_{j2s},\cdots,y_{jr_js})^T$
is a vector in $\mathbf{R}^{r_j}$.

Finally, since $E_{\alpha,\alpha}
(\lambda_jt^\alpha)>0$ for all $t\geq 0$, $j=1,2,\cdots,$ we then show our proof by using the Reductio and Absurdum.

$(a)$ Necessity. If $p\geq r=\max\{r_j\}$ and $rank ~G_j^s<r_j \mbox{ for some } j=1,2,\cdots$ and $s=1,2,\cdots,p,$ there exists a nonzero
element $\tilde{y} \in \left( L^2( \omega)\right)^n$ with
$\tilde{y}_{js}=\left(\tilde{y}_{j1s},\tilde{y}_{j2s},\cdots, \tilde{y}_{jr_js}\right)^T\in \mathbf{R}^{r_j}$
 such that
$
G_{j}^s\tilde{y}_{js}=\mathbf{0}_p.
$
Then we can find a nonzero vector  $\tilde{y} \in \left( L^2( \omega)\right)^n$ satisfying
$
\sum\limits_{j=1}^{\infty}{E_{\alpha,\alpha}(\lambda_jt^{\alpha})}\sum\limits_{s=1}^n G_j^s\tilde{y}_{js}=\mathbf{0}_p.
$
This means that the sensors $(D_i,f_i)_{1\leq i\leq p}$ are not $\omega-$strategic.

$(b)$ Sufficiency.  On the contrary, if the sensors $(D_i,f_i)_{1\leq i\leq p}$ are not $\omega-$strategic, i.e.,
$
\overline{Im\left( p_{\omega}\nabla K^*\right)}\neq \left(L^2(\omega)\right)^n.
$
Then there exists a nonzero element  $y_{j^*s}\in \mathbf{R}^{r_j}$ such that
$
\sum\limits_{s=1}^n G_{j^*}^sy_{j^*s}=\mathbf{0}_p.
$
This allows us to complete the first conclusion of the theorem.

In particular,  when $n=s=1,$  similar to the argument in $(a)$, if $p\geq r=\max\{r_j\}$ and $rank ~G_j^1<r_j \mbox{ for some } j=1,2,\cdots$,
there  exists a nonzero vector  $\tilde{y} \in \left( L^2( \omega)\right)^n$ satisfying
$
\sum\limits_{j=1}^{\infty}{E_{\alpha,\alpha}(\lambda_jt^{\alpha})} G_j^1\tilde{y}_{js}=\mathbf{0}_p.
$
Then the sensors $(D_i,f_i)_{1\leq i\leq p}$ are not $\omega-$strategic.

Moreover, if the sensors $(D_i,f_i)_{1\leq i\leq p}$ are not $\omega-$strategic, there exists a nonzero element  $y\in \left( L^2( \omega)\right)^n$ satisfying
$
 G_j^1y_{j1}=\mathbf{0}_p.
$
Then if $p\geq r=\max\{r_j\}$, it is sufficient to see that
$
rank ~G_{j}^1<r_{j }
$
for all $j=1,2,\cdots$.
The proof is complete.

\section{An approach for the regional gradient  reconstruction}
\label{section4}

This section is focused on an approach, which allows us to reconstruct  the initial gradient vector of the system $(\ref{problem})$ in $\omega.$ The method used here
is Hilbert uniqueness method (HUMs) introduced by \cite{Lions}, which can be considered as  an extension of those given in \cite{rgb2}.

Let  $G$ be the set given by
\begin{eqnarray*}
G=\left\{
\begin{array}{l}
g\in \left(L^2(\Omega)\right)^n: g=\mathbf{0}_n\mbox{ in } \Omega\backslash\omega \mbox{ and there exists} \\{\kern 37pt}\mbox{ a unique }
 \tilde{g}\in H^1_0(\Omega)
\mbox{ such that } \nabla \tilde{g}=g
\end{array}\right\}.
\end{eqnarray*}
For any $g^*\in G$, there exists a function $\tilde{g}^*\in H^1_0(\Omega)$ satisfying $\tilde{g}^*=\nabla^*p_\omega^*g^*$. Consider the following system
\begin{equation}\label{Gfunction}
\left\{\begin{array}{l}
{}_0D^{\alpha}_t\varphi(t)=A\varphi(t), ~~ t\in [0,b],
\\ \lim\limits_{t\to 0^+} {}_0I^{1-\alpha}_t\varphi(t)=\tilde{g}^*,
\end{array}\right.
\end{equation}
which admits a unique solution $\varphi\in L^2(0,b; H^1_0(\Omega))\cap C([0,b]\times \Omega)$ given by
$
\varphi(t)=S_\alpha(t)\tilde{g}^*=S_\alpha(t)\nabla^*p_\omega^*g^* .
$
Then we consider the  semi-norm on $G$
\begin{equation}\label{Gnorm}
g^*\in G \to \|g^*\|_G^2=\int_0^b{\| C\varphi (b-t)\|^2}dt
\end{equation}
and we can get the following result.
\begin{lem}\label{lemma4.1}
If the system $(\ref{problem})$ is  regionally gradient observable in $\omega$ at time $b$, then $(\ref{Gnorm})$ defines a norm on $G$.
\end{lem}
$\textbf{Proof.}$
 If the system $(\ref{problem})$ is  regionally gradient observable, by Definition $\ref{rgodef}$,  one has
\begin{eqnarray*}
Ker  \left(K(t)\nabla^* p^*_\omega\right)&=&Ker  \left(CS_\alpha(t) \nabla^*p_\omega^*\right)\\&=& \{\mathbf{0}_n\}\subseteq  \left(L^2(\omega)\right)^n.\end{eqnarray*}
Moreover,  for any $g^*\in G,$ since
\begin{eqnarray*}
  &~&\|g^*\|_G=0
\Leftrightarrow CS_\alpha(b-t)\nabla^* p^*_\omega g^*=0,~~\forall t\in(0,b),\end{eqnarray*}
which gives   $g^*= \mathbf{0}_n \in \left(L^2(\omega)\right)^n, $ it then follows that $(\ref{Gnorm})$ defines  a norm of $G$ and the proof is complete.

For $g^*\in G$, consider the operator $\Lambda: G\to G^*$ defined by
\begin{eqnarray}\label{4.5}
\Lambda g^*=p_\omega  \nabla \psi(b),
\end{eqnarray}
where $ \psi(t)$ solves
the following system $(t\in [0,b])$
\begin{equation}
\left\{\begin{array}{l}
_0D^{\alpha}_{t}\psi(t)=A^*\psi(t)+C^*C\varphi(b-t),
\\ \lim\limits_{t\to 0^+}{}_0I^{1-\alpha}_{t}\psi(t)=0
\end{array}\right.
\end{equation}
 controlled by the solution of the system $(\ref{Gfunction})$.
We then conclude that the regional gradient reconstruction  problem is equivalent to solving the equation $(\ref{4.5})$.

\begin{thm}\label{Theorem4.1}
If $(\ref{problem})$ is  regionally gradient observable in $\omega$ at time $b$, then $(\ref{4.5})$ has a unique
solution $g^*\in G$ and the initial gradient $\nabla y_0$ in subregion $\omega$
is equivalent to $g^*$.
\end{thm}
$\textbf{Proof.}$
By Lemma $\ref{lemma4.1}$,  we see that  $\|\cdot\|_G$ is a norm of the space $G$ provided that the system $(\ref{problem})$ is
regionally gradient observable in $\omega$ at time $b$.
Let the completion of $G$ with respect to the norm $\|\cdot\|_G$ again by $G$. By the  Theorem 1.1 in \cite{Lions2},
to obtain the existence of the unique solution $g^*\in G$ of problem $(\ref{4.5})$, we only need to show that $\Lambda$ is coercive from $G$ to $G^*,$
i.e., there exists a constant $\mu>0$ such that
\begin{eqnarray}
(\Lambda g,g)_{\left(L^2(\Omega)\right)^n}\geq \mu\|g\|_G^2,~~\forall g\in G.
\end{eqnarray}
Indeed, for any $g^*\in G,$ we have
\begin{eqnarray*}
&~&(\Lambda g^*,g^*)_{\left(L^2(\Omega)\right)^n}=(p_\omega  \nabla \psi(b) ,g^*)_{\left(L^2(\Omega)\right)^n} \\
&=&\left( \nabla \int_0^bS^*_\alpha(b-s)C^*C\varphi(b-s)ds,p_\omega^* g^*\right)_{\left(L^2(\Omega)\right)^n}\\
&=&\int_0^b\left( C\varphi(b-s),CS_\alpha(b-s)\nabla^*p_\omega^*g^*\right)ds\\
&=&\|g^*\|_G^{2},
\end{eqnarray*}
Then $\Lambda$ is coercive  and $(\ref{4.5})$ has a unique solution,
which is also the initial gradient to be estimated in the subregion $\omega$ at time $b$.  The proof is complete.

\begin{rem}
Note that if the Riemann-Liouville fractional derivative in system $(1)$
is replaced by  a Caputo fractional derivative,  its unique mild solution will be given by (\cite{Japan})
\begin{eqnarray}
y(t)=\eta_\alpha(t)y_0,~~\eta_\alpha(t)\neq S_\alpha(t),~t\in [0,b].
\end{eqnarray}
We see  that $z(t)=C\eta_\alpha(t)y_0=K(t)y_0$ and
\begin{equation}
K^*:~\left\{\begin{array}{l}L^2(0,b;\mathbf{R}^p)\to Y,
\\~ z\to \int_0^b{\eta^*_\alpha(s)C^*z(s)}ds.
\end{array}\right.
\end{equation}
Then the Lemma 7 fails. New lemmas similar to Lemma 4 and Lemma 7 are of great interest.
Besides, this challenge is also our interest now and we shall try our best
to study it in our forthcoming papers.
\end{rem}
\section{Applications}
Let $\Omega_2=[0,1]\times[0,1]$ and $\omega_2\subseteq \Omega_2$. In this section, let us consider the following system
\begin{eqnarray}\label{example}
\left\{\begin{array}{l}
_0D^{\alpha}_{t}y(x,t)=\triangle y(x,t)~\mbox{ in }\Omega_2 \times [0,b],
\\ y(\eta, t)= 0~\mbox{ on }\partial\Omega_2 \times [0,b],\\
\lim\limits_{t\to 0^+}{}_0I^{1-\alpha}_{t}y(x,t)=y_0(x)~\mbox{ in }\Omega_2,
\end{array}\right.
\end{eqnarray}
where $\triangle $ is the elliptic operator given by
$
\triangle=\frac{\partial^2 }{\partial  x_1^2}+\frac{\partial^2  }{\partial x_2^2 }.
$
For any $i$$(i=1,2)$, let the gradient of initial vector be
\begin{eqnarray}
\frac{\partial y_0}{\partial x_i}=\left\{\begin{array}{l}
\frac{\partial y^1_0}{\partial x_i}~\mbox{ on }\omega_2 ,
\\ \frac{\partial y^2_0}{\partial x_i}~\mbox{ on }\Omega_2\backslash \omega_2.
\end{array}\right.
\end{eqnarray}
Then our aim here  is to present an approach to reconstruct the regional gradient vector:
$\left(\frac{\partial y^1_0}{\partial x_1},\frac{\partial y^1_0}{\partial x_2}\right),$
where the sensors may be zone, pointwise  or filament ones.\\\\
\textbf{Case 1. Zone sensors}\\
Suppose that $D=[d_1,d_2]\times [d_3,d_4]\subseteq \Omega_2 $ and the output functions are
\begin{eqnarray*}
z(t)=Cy(x,t)=\int_{D}{f(x)y(x,t)}dx,f\in  L^2(D),x\in \Omega_2,
\end{eqnarray*}
where the system is observed by one sensor $(D,f)$ and $C$ is bounded.
Moreover, we get that the eigenvalue, corresponding  eigenvector of $\triangle$ and the semigroup generated by $\triangle$ are
 $\lambda_{mn}=-(m^2+n^2)\pi^2$, $\xi_{mn}=2\sin (m\pi x_1)\sin (n\pi x_2)$ and $~y\in L^2(\Omega_2),$\\
$\Phi(t)y(x)=\sum\limits_{m,n=1}^{\infty}{\exp(\lambda_{mn}t)(y,\xi_{mn})}\xi_{mn}(x),$
respectively. Then the multiplicity of the eigenvalues is one and
\begin{eqnarray*}\label{K2}
\begin{array}{l}
K(t)y_0=CS_\alpha(t)y_0=\sum\limits_{m,n=1}^{\infty}\\\frac{E_{\alpha,\alpha}(\lambda_{mn}t^{\alpha})
(y_0,\xi_{mn})}{t^{\alpha-1}}\int_{d_3}^{d_4}{\int_{d_1}^{d_2}f(x_1,x_2)\xi_{mn}(x_1,x_2)dx_1}dx_2.
\end{array}
\end{eqnarray*}
Let $f(x_1,x_2)=\sin(\sqrt{2}\pi x_1)\sin(\sqrt{2}\pi x_2)$. We see that $(\ref{example})$  is not gradient observable on $\Omega_2$. However,
\begin{prop}\label{prop12}
The sensor $(D,f)$ is gradient strategic in $\omega_2\subseteq \Omega_2$ if and only if
\[\begin{array}{l}\beta_{1m} \left(p_{1\omega_2}y_1,\xi_{mn}\right)+\beta_{2n}\left(p_{1\omega_2}y_2,\xi_{mn}\right)= 0,~ \forall m,n=1,2,\cdots  \\
\Rightarrow (y_1,y_2)=(0,0),\end{array}\]
where

$\beta_{1m}=m\pi \int_{d_1}^{d_2}\int_{d_3}^{d_4}{f(x_1,x_2)\cos( m \pi x_1)\sin( n\pi x_2)}dx_1dx_2$ and
$\beta_{2n}=n\pi\int_{d_1}^{d_2}\int_{d_3}^{d_4}{f(x_1,x_2)\sin( m\pi x_1)\cos( n\pi x_2)}dx_1dx_2.$
\end{prop}
\textbf{Proof.}
According to the argument above, we have $p=1,~r_j=1,~n=2$. It then follows that
$G_j^1=2 \left[\beta_{1m}\right]_{1 \times 1} $ and $ G_j^2=2 \left[\beta_{2n}\right]_{1 \times 1}.$
Let $y_{j1}=\left(p_{1\omega_2}y_1,\xi_{mn}\right)$ and $y_{j2}=\left(p_{1\omega_2}y_2,\xi_{mn}\right).$
By Theorem $\ref{theorem3.1}$,
then the necessary and sufficient condition for the sensor $(D,f)$ to be gradient strategic in $\omega_2\subseteq \Omega_2$
is that
$\begin{array}{l}\beta_{1m} y_{j1}+\beta_{2n}y_{j2}= 0 \mbox{ for all } m,n=1,2,\cdots  \\
\Rightarrow (y_1,y_2)=(0,0). \end{array}$\\The proof is complete.

Let  $G_2$ be a set defined by
\begin{equation}\label{G2}
G_2=\left\{\begin{array}{l}g\in \left(L^2(\Omega_2)\right)^2: g=0\mbox{ in } \Omega_2\backslash\omega_2 \mbox{ and }
\\ {\kern 30pt}\mbox{there exists a unique }
 \tilde{g}\in H^1_0(\Omega_2)
\\{\kern 30pt} \mbox{ such that } \nabla \tilde{g}=g
\end{array}\right\}
\end{equation}
and for any $g^*\in G_2$, we see that \\
$\|g^*\|_{G_2}^2=\int_0^b{\left[\left(S_\alpha(b-t) \sum\limits_{s=1}^2\frac{\partial(p_\omega^*g^*)}{\partial x_s},f\right)_{L^2(D)}\right]^2}dt$
\\defines a norm on $G_2$ provided that $(\ref{example})$ is regionally gradient observable in $\omega_2$ at time $b$.
Consider the system
\begin{equation}
\left\{\begin{array}{l}
_0D^{\alpha}_{t}\psi(x,t)=\frac{\partial^2}{\partial x^2}\psi(x,t)+p_D f(x)\\
{\kern 6pt}\times \left(S_\alpha(b-t) \sum\limits_{s=1}^2\frac{\partial(p_\omega^*g^*)}{\partial x_s},f\right)_{L^2(D)}\mbox{ in }\Omega_2 \times [0,b],\\
\psi(\eta, t)=0~\mbox{ on }\partial\Omega_2 \times [0,b],\\
\lim\limits_{t\to 0^+}{}_0I^{1-\alpha}_{t}\psi(x,t)=0\mbox{ in }\Omega_2.
\end{array}\right.
\end{equation}
By Theorem $\ref{Theorem4.1}$, the  equation $\Lambda: g^* \to p_\omega  \nabla \psi(b)$
has a unique solution in $G_2$, which is also the initial gradient $\left(\frac{\partial y_0}{\partial x_1},\frac{\partial y_0}{\partial x_2}\right)$ in $\omega_2.$\\\\
\textbf{Case 2. Pointwise sensors}\\
In this part, we consider the problem $(\ref{example})$ with the following unbounded output function
\begin{eqnarray}
z(t)=Cy(x,t)=y(\sigma,t),~~\sigma\in  \Omega.
\end{eqnarray}

Let $\sigma=(\sigma_1,\sigma_2)\in \Omega_2$ be the location of the sensor and let $d_1=d_2=\sigma_1$, $d_3=d_4=\sigma_2$ in Eq $(\ref{K2})$, then one has
\begin{eqnarray}
K(t)y_0=\sum\limits_{m,n=1}^{\infty}t^{\alpha-1}E_{\alpha,\alpha}(\lambda_{mn}t^{\alpha})(y_0,\xi_{mn})\xi_{mn}(\sigma).
\end{eqnarray}
Since $|\xi_{mn}|\leq 2~\mbox{for }x\in [0,1]\times[0,1]$, $E_{\alpha,\alpha}(\lambda_{mn}t^{\alpha})$ is continuous and
$\left|E_{\alpha,\alpha}(\lambda_{mn}t^{\alpha})\right|\leq \frac{C}{1+|\lambda_{mn}|t^\alpha}~(C>0)$ (\cite{Podlubny}),
we get that the assumption $(A_1)$ is satisfied.
Further, for any $z\in L^2(\Omega),$ one has
\begin{eqnarray*}
\begin{array}{l}
K^*z(t)=\sum\limits_{m,n=1}^{\infty}\int_0^b{ E_{\alpha,\alpha}(\lambda_{mn}\tau^{\alpha}) (C^*z(\tau),\xi_{mn})}d\tau\xi_{mn}(x).
\end{array}\end{eqnarray*}
Then the assumption $(A_2)$ holds.
By Theorem $\ref{theorem3.1}$,
similar to Proposition $\ref{prop12}$, let $d_1=d_2=\sigma_1$ and $d_3=d_4=\sigma_2$, we see that,
\begin{prop}
There exists a subregion $\omega_2\subseteq \Omega_2$ such that the sensor $(\sigma,\delta_{\sigma})$ is gradient $\omega_2-$strategic if and only if
$m\pi\cos( m \pi \sigma_1)\sin( n\pi\sigma_2)\left(p_{1\omega_2}y_1,\xi_{mn}\right)
+ n\pi\sin( m \pi \sigma_1)\cos( n\pi\sigma_2)\left(p_{1\omega_2}y_2,\xi_{mn}\right)= 0, ~~\forall m,n=1,2,\cdots$
can imply $(y_1,y_2)=(0,0)$.
\end{prop}
Further, for any $g^*\in G_2$, by Lemma $\ref{lemma4.1},$ if $(\ref{example})$ is  regionally gradient observable, then
\[
 \|g^*\|_{G_2}^2=\int_0^b{\left[\left(S_\alpha(b-t) \sum\limits_{s=1}^2\frac{\partial(p_\omega^*g^*)}{\partial x_s}\right)(\sigma)\right]^2}dt
\]
defines a norm on $G_2$. Consider the following system
\begin{equation}
\left\{\begin{array}{l}
_0D^{\alpha}_{t}\psi(x,t)=\frac{\partial^2}{\partial x^2}\psi(x,t)+\delta(x-\sigma)\\
{\kern 18pt}\times \left(S_\alpha(b-t) \sum\limits_{s=1}^2\frac{\partial(p_\omega^*g^*)}{\partial x_s}\right),~(x,t)\in \Omega_2 \times [0,b],\\
\psi(\eta , t)=0,~~(\eta,t)\in \partial\Omega_2 \times [0,b],\\
\lim\limits_{t\to 0^+}{}_0I^{1-\alpha}_{t}\psi(x,t)=0, ~~x\in \Omega_2.
\end{array}\right.
\end{equation}
It follows from Theorem $\ref{Theorem4.1}$ that the  equation $\Lambda: g^* \to p_\omega  \nabla \psi(b)$
has a unique solution in $G_2$, which is also the initial gradient $\left(\frac{\partial y_0}{\partial x_1},\frac{\partial y_0}{\partial x_2}\right)$ on $\omega_2.$\\\\
\textbf{Case 3. Filament  sensors}\\
Consider the case where the observer $(F, \delta_F)$ is located on  the curve
 $F=[\tau_1,\tau_2]\times \{\sigma\}\subseteq \Omega_2 $ and the output functions  are
$
z(t)=\int_{\tau_1}^{\tau_2}{\delta_F(x_1,\sigma)y(x_1,\sigma,t)}dx_1.
$
For example, let $\delta_F(x_1,x_2)=\sin(\sqrt{2}\pi x_1)\sin(\pi x_2)$. Then
the example $(\ref{example})$  is not gradient observable in $\Omega_2$  at time $b$.
 By Theorem $\ref{theorem3.1}$,
let $d_1=\tau_1,~d_2=\tau_2$ and $d_3=d_4=\sigma$ in Proposition $\ref{prop12}$, we get the following results.
\begin{prop}
The sensor $(F,\delta_F)$ is gradient strategic in a subregion $\omega_2\subseteq \Omega_2$ if and only if
\[\begin{array}{l}\gamma_{1m} \left(p_{1\omega_2}y_1,\xi_{mn}\right)+\gamma_{2n}\left(p_{1\omega_2}y_2,\xi_{mn}\right)= 0,\forall m,n=1,2,\cdots \\ \Rightarrow (y_1,y_2)=(0,0),\end{array}\]
where
$\gamma_{1m}=m\pi\int_{\tau_1}^{\tau_2}{\delta_F(x_1,\sigma)\cos( m \pi x_1)\sin( n\pi \sigma)}dx_1$ and
$\gamma_{2n}=n\pi \int_{\tau_1}^{\tau_2}{\delta_F(x_1,\sigma)\sin( m \pi x_1)\cos( n\pi \sigma)}dx_1,$
for all $m,n=1,2,\cdots$.
\end{prop}

Let  $G_2$ be defined by $(\ref{G2})$ and for any $g^*\in G_2$, consider
\begin{equation}
\left\{\begin{array}{l}
_0D^{\alpha}_{t}\psi(x,t)=\frac{\partial^2}{\partial x^2}\psi(x,t)+p_F \delta_F(x)\\
{\kern 10pt}\times \left(S_\alpha(b-t) \sum\limits_{s=1}^2\frac{\partial(p_\omega^*g^*)}{\partial x_s},\delta_F\right)_{L^2(F)} \mbox{ in }\Omega_2 \times [0,b],\\
\psi(\eta, t)=0~\mbox{ on }\partial\Omega_2 \times [0,b],\\
\lim\limits_{t\to 0^+}{}_0I^{1-\alpha}_{t}\psi(x,t)=0\mbox{ in }\Omega_2,
\end{array}\right.
\end{equation}
where $
 \|g^*\|_{G_2}^2=\int_0^b{\left[\left(S_\alpha(b-t) \sum\limits_{s=1}^2\frac{\partial(p_F^*g^*)}{\partial x_s},\delta_F\right)_{L^2(F)}\right]^2}dt
$
is a norm on $G_2$ and by Theorem $\ref{Theorem4.1}$, the  equation $\Lambda: g^* \to p_\omega  \nabla \psi(b)$
has a unique solution in $G_2$ and $g^*=\left(\frac{\partial y_0}{\partial x_1},\frac{\partial y_0}{\partial x_2}\right)$ on $\omega_2.$

\section{CONCLUSIONS}
In this paper, we investigate the regional gradient observability problem for the
time fractional diffusion system with  Riemann-Liouville fractional derivatie, which is motivated
by many real world applications where the objective is to obtain useful information on the state
gradient in a given subregion of the whole domain. We hope that the results here could provided
some insights into the control theoretical analysis of fractional order systems.
Moreover, the results presented
here can also be extended to complex fractional order DPSs and  various
open questions are still under consideration. For example, the problem of state gradient control
of fractional order DPSs, regional observability of fractional order
system with mobile sensors as well as the regional sensing configuration
are of great interest. For more
information on the potential topics related to fractional DPSs, we
refer the readers to \cite{ge2015JAS} and the references therein

\end{document}